\magnification\magstep1
\baselineskip 16pt
\def\pfbox
  {\hbox{\hskip 3pt\lower2pt\vbox{\hrule
  \hbox to 5pt{\vrule height 7pt\hfill\vrule}
  \hrule}}\hskip3pt}
\def\proof{\noindent{\bf Proof.}\enspace}
\def\bib{\par\smallskip\noindent\hangindent 25pt}

\centerline{\bf Johann Faulhaber and Sums of Powers}
\centerline{\sl Donald E. Knuth}
\centerline{\sl Computer Science Department, Stanford University}

\vfill

{\narrower\smallskip\noindent
{\bf Abstract.}\enspace
Early 17th-century mathematical publications of Johann Faulhaber
contain some remarkable theorems, such as the fact that the $r$-fold
summation of $1^m,2^m,\ldots,n^m$ is a polynomial in $n(n+r)$ when $m$
is a positive odd number. The present paper explores a
computation-based approach by which Faulhaber may well have discovered
such results, and solves a 360-year-old riddle that Faulhaber
presented to his readers. It also shows that 
similar results hold when we express the sums in terms of central
factorial powers instead of ordinary powers. Faulhaber's coefficients
can moreover be generalized to factorial powers of
noninteger exponents, obtaining asymptotic series for
$1^{\alpha}+2^{\alpha}+\cdots+n^{\alpha}$
in powers of $n^{-1}(n+1)^{-1}$.
\smallskip}

\vfill\vfill\eject

\bigskip\noindent
Johann Faulhaber of Ulm (1580--1635),
 founder of a school for engineers early in
the 17th century, loved numbers. His passion for arithmetic and
algebra led him to devote a considerable portion of his life to the
computation of formulas for the sums of powers, significantly
extending all previously known results. 
He may well have carried out
 more computing than anybody else in Europe during the
first half of the 17th~century. His greatest mathematical
achievements appear in a booklet entitled {\sl Academia Algebr{\ae}\/}
(written in German in spite of its latin title), published in
Augsburg, 1631~[2].
Here we find, for example, the following formulas for sums of odd
powers:
$$\eqalign{%
1^1+2^1+\cdots+n^1&=N\,,\qquad N=(n^2+n)/2\,;\cr
1^3+2^3+\cdots+n^3&=N^2\,;\cr
1^5+2^5+\cdots +n^5&=(4N^3-N^2)/3\,;\cr
1^7+2^7+\cdots +n^7&=(12N^4-8N^3+2N^2)/6\,;\cr
1^9+2^9+\cdots+n^9&=(16N^5-20N^4+12N^3-3N^2)/5\,;\cr
1^{11}+2^{11}+\cdots+n^{11}&=(32N^6-64N^5+68N^4-40N^3+5N^2)/6\,;\cr
1^{13}+2^{13}+\cdots+n^{13}
&=(960N^7-2800N^6+4592N^5-4720N^4+2764N^3\cr
&\qquad\null-691N^2)/105\,;\cr
1^{15}+2^{15}+\cdots+n^{15}
&=(192N^8-768N^7+1792N^6-2816N^5+2872N^4\cr
&\qquad\null-1680N^3+420N^2)/12\,;\cr
1^{17}+2^{17}+\cdots+n^{17}
&=(1280N^9-6720N^8+21120N^7-46880N^6+72912N^5\cr
&\qquad\null-74220N^4+43404N^3-
10851N^2)/45\,.\cr}$$
Other mathematicians had studied $\Sigma n^1,\Sigma n^2,\ldots,\,\Sigma
n^7$  and he had previously gotten as far as $\Sigma n^{12}$;
 but the sums had always
previously been expressed as polynomials in~$n$, not~$N$.

Faulhaber begins his book by simply stating these novel formulas and
proceeding to expand them into the corresponding polynomials in~$n$.
Then he verifies the results when $n=4$, $N=10$. But he gives no clues
about how he derived the expressions; he states only that the leading
coefficient in $\Sigma n^{2m-1}$ will be $2^{m-1}\!/m$, and that the trailing
coefficients will have the form $4\alpha_mN^3-\alpha_mN^2$ when
$m\geq 3$. 

Faulhaber believed that similar polynomials in~$N$, with alternating
signs, would continue to exist for all~$m$, but he may not really have
known how to prove such a theorem. In his day, mathematics was treated
like all other sciences; it was sufficient to present a large body of
evidence for an observed phenomenon. A~rigorous proof of Faulhaber's
assertion was first
published by Jacobi in 1834~[6].
A.~W.~F. Edwards showed recently how to obtain the coefficients by
matrix inversion~[1],
based on another proof given by L.~Tits in 1923~[8].
But none of these proofs use methods that are very close to those
known in 1631.

Faulhaber went on to consider sums of sums. Let us write $\Sigma^rn^m$
for the $r$-fold summation of $m$th powers from~1 to~$n$; thus,
$$\Sigma^0n^m=n^m\,;\qquad
\Sigma^{r+1}n^m=\Sigma^r1^m+\Sigma^r2^m+\cdots +\Sigma^rn^m\,.$$
He discovered that $\Sigma^rn^{2m}$ can be written as a polynomial in
the quantity
$$N_r=(n^2+rn)/2\,,$$
times $\Sigma^rn^2$. For example, he gave the formulas
$$\vcenter{\halign{$#$\hfil\ &$#$\hfil\cr
\Sigma^2n^4&=(4N_2-1)\,\Sigma^2n^2\!/5\,;\cr
\Sigma^3n^4&=(4N_3-1)\,\Sigma^3n^2\!/7\,;\cr
\Sigma^4n^4&=(6N_4-1)\,\,\Sigma^4n^2\!/14\,;\cr
\Sigma^6n^4&=(4N_6+1)\,\Sigma^6n^2\!/15\,;\cr
\Sigma^2n^6&=(6N_2^2-5N_2+1)\,\Sigma^2n^2\!/7\,;\cr
\Sigma^3n^6&=(10N_3^2-10N_3+1)\,\Sigma^3n^2\!/21\,;\cr
\Sigma^4n^6&=(4N_4^2-4N_4-1)\,\Sigma^4n^2\!/14\,;\cr
\Sigma^2n^8&=(16N_2^3-28N_2^2+18N_2-3)\,\Sigma^2n^2\!/15\,.\cr}}$$
He also gave similar formulas for odd exponents, factoring out
$\Sigma^rn^1$ instead of~$\Sigma^rn^2$:
$$\eqalign{\Sigma^2n^5&=(8N_2^2-2N_2-1)\Sigma^2n^1\!/14\,;\cr
\Sigma^2n^7&=(40N_2^3-40N_2^2+6N_2+6)\Sigma^2n^1\!/60\,.\cr}$$
And he claimed that, in general, $\Sigma^rn^m$ can be expressed as a
polynomial in~$N_r$ times either $\Sigma^rn^2$ or $\Sigma^rn^1$,
depending on whether $m$ is even or odd.

Faulhaber had probably verified this remarkable theorem in many cases
including $\Sigma^{11}n^6$, because he exhibited a polynomial in~$n$
for $\Sigma^{11}n^6$ that would have been quite difficult to obtain by
repeated summation. His polynomial, which has the form
$${6n^{17}+561n^{16}+\cdots+1021675563656n^5+\cdots-96598656000n\over
2964061900800}\,,$$
turns out to be absolutely correct,
according to calculations with a modern computer. (The denominator is
$17!/120$. One cannot help
thinking that nobody has ever checked these numbers since Faulhaber
himself wrote them down, until today.)

Did he, however, know how to prove his claim, in the sense that
20th~century mathematicians would regard his argument as conclusive?
He may in fact have known how to do so, because there is an extremely
simple way to verify the result using only methods that he would have
found natural.

\bigskip\noindent
{\bf Reflective functions.}\enspace
Let us begin by studying an elementary property of integer functions.
We will say that the function $f(x)$ is $r$-{\it reflective\/} if
$$f(x)=f(y)\quad{\rm whenever}\quad x+y+r=0\,;$$
and it is {\it anti\/}-$r$-{\it reflective\/} if
$$f(x)=-f(y)\quad{\rm whenever}\quad x+y+r=0\,.$$
The values of $x$, $y$, $r$ will be assumed to be integers for
simplicity. 
When $r=0$, reflective functions are
even, and anti-reflective functions are odd.
Notice that $r$-reflective functions are closed under
addition and multiplication; the product of two anti-$r$-reflective
functions is $r$-reflective. 

Given a function $f$, we define its backward difference $\nabla f$ in
the usual way:
$$\nabla f(x)=f(x)-f(x-1)\,.$$
It is now easy to verify a simple basic fact.

\proclaim
Lemma 1. If $f$ is $r$-reflective then $\nabla f$ is
anti-$(r-1)$-reflective. If $f$ is anti-$r$-reflective then $\nabla f$
is $(r-1)$-reflective.

\proof
If $x+y+(r-1)=0$ then $x+(y-1)+r=0$ and $(x-1)+y+r=0$. Thus $f(x)=\pm
f(y-1)$ and $f(x-1)=\pm f(y)$ when $f$ is $r$-reflective or
anti-$r$-reflective. 
\pfbox

\bigskip
Faulhaber almost certainly knew this lemma, because [2, folio D.iii recto]
presents a table of $n^8,\nabla n^8,\ldots,\nabla^8n^8$ in which the
reflection phenomenon is clearly apparent. He states that he has
constructed ``grosse Tafeln,'' but that this example should be ``alles
gnugsam vor Augen sehen und auf h\"ohere quantiteten [exponents]
continuiren k\"onde.''

The converse of Lemma~1 is also true, if we are careful. Let us
define $\Sigma$ as an inverse to the $\nabla$~operator:
$$\Sigma f(n)=\cases{C+f(1)+\cdots+f(n)\,,&if $n\geq 0$;\cr
\noalign{\smallskip}
C-f(0)-\cdots-f(n+1)\,,&if $n<0$.\cr}$$
Here $C$ is an unspecified constant, which we will choose later;
whatever its value, we have
$$\nabla\Sigma f(n)=\Sigma f(n)-\Sigma f(n-1)=f(n)$$
for all $n$.

\proclaim
Lemma 2. If $f$ is $r$-reflective, there is a unique $C$ such that
$\Sigma f$ is anti-$(r+1)$-reflective. If $f$ is anti-$r$-reflective,
then $\Sigma f$ is $(r+1)$-reflective for all~$C$.

\proof
If $r$ is odd, $\Sigma f$ can be anti-$(r+1)$-reflective only if $C$
is chosen so that we have  $\Sigma f\bigl(-(r+1)/2\bigr)=0$. If $r$ is even,
$\Sigma f$ can be anti-$(r+1)$-reflective only if $\Sigma
f(-r/2)=-\Sigma f(-r/2-1)=-\bigl(\Sigma f(-r/2)-f(-r/2)\bigr)$; 
i.e., $\Sigma f(-r/2)={1\over 2}f(-r/2)$.

Once we have found $x$ and~$y$ such that $x+y+r+1=0$ and $\Sigma
f(x)=-\Sigma f(y)$, it is easy to see that we will also have $\Sigma
f(x-1)=-\Sigma f(y+1)$, if $f$ is $r$-reflective, since $\Sigma
f(x)-\Sigma f(x-1)=f(x)=f(y+1)=\Sigma f(y+1)-\Sigma f(y)$.

Suppose on the other hand that $f$ is anti-$r$-reflective. If $r$ is
odd, clearly $\Sigma f(x)=\Sigma f(y)$ if $x=y=-(r+1)/2$. If $r$ is
even, then $f(-r/2)=0$; so $\Sigma f(x)=\Sigma f(y)$ when
$x=-r/2$ and $y=-r/2-1$. Once we have found $x$ and~$y$ such that
$x+y+r+1=0$ and $\Sigma f(x)=\Sigma f(y)$, it is easy to verify as
above that $\Sigma f(x-1)=\Sigma f(y+1)$.\quad\pfbox

\proclaim
Lemma 3. If $f$ is any even function with $f(0)=0$, the $r$-fold
repeated sum $\Sigma^rf$ is $r$-reflective for all even~$r$ and
anti-$r$-reflective for all odd~$r$, if we choose the constant $C=0$
in each summation. If $f$ is any odd function, the $r$-fold repeated
sum $\Sigma^rf$ is $r$-reflective for all odd~$r$ and
anti-$r$-reflective for all even~$r$, if we choose the constant $C=0$
in each summation.

\proof
Note that $f(0)=0$ if $f$ is odd. If $f(0)=0$ and if we
always choose $C=0$, it is easy to verify by induction on~$r$ that
$\Sigma^rf(x)=0$ for $-r\leq x\leq 0$. Therefore the choice $C=0$
always agrees with the unique choice stipulated in the proof of
Lemma~2, whenever a specific value of~$C$ is necessary in that
lemma.\quad\pfbox

\bigskip
When $m$ is a positive integer, the function $f(x)=x^m$ obviously
satisfies the condition of Lemma~3. Therefore we have proved that each
function $\Sigma^rn^m$ is either $r$-reflective or
anti-$r$-reflective, for all $r>0$ and $m>0$. And Faulhaber presumably
knew this too. His theorem can now be proved if we supply one small
additional fact, specializing from arbitrary functions to polynomials:

\proclaim
Lemma 4. A polynomial $f(x)$ is $r$-reflective if and only if it can
be written as a polynomial in $x(x+r)$; it is anti-$r$-reflective if
and only if it can be written as $(x+r/2)$ times a polynomial in
$x(x+r)$.

\proof
The second statement follows from the first, because we  have already
observed that an anti-$r$-reflective function must have $f(-r/2)=0$
and because the function $x+r/2$ is obviously anti-$r$-reflective.
Furthermore, any polynomial in $x(x+r)$ is $r$-reflective, because
$x(x+r)=y(y+r)$ when $x+y+r=0$. Conversely, if $f(x)$ is
$r$-reflective we have $f(x-r/2)=f(-x-r/2)$, so $g(x)=f(x-r/2)$ is an
even function of~$x$; hence $g(x)=h(x^2)$ for some polynomial~$h$.
Then $f(x)=g(x+r/2)=h\bigl(x(x+r)+r^2/4\bigr)$ is a polynomial in
$x(x+r)$.\quad \pfbox

\proclaim
Theorem {\rm (Faulhaber)}. 
There exist polynomials $g_{r,m}$ for all positive integers~$r$
and~$m$ such that
$$\Sigma^rn^{2m-1}=g_{r,2m+1}\bigl(n(n+r)\bigr)\Sigma^rn^1\,,\qquad
\Sigma^rn^{2m}=g_{r,2m}\bigl(n(n+r)\bigr)\Sigma^rn^2\,.$$

\proof
Lemma 3 tells us that $\Sigma^rn^m$ is $r$-reflective if $m+r$ is even
and anti-$r$-reflective if $m+r$ is odd. 

Note that $\Sigma^rn^1={n+r\choose r+1}$. Therefore a polynomial
in~$n$ is a multiple of $\Sigma^rn^1$ if and only if it vanishes at
$-r,\ldots,-1,0$. We have shown in the proof of Lemma~3 that
$\Sigma^rn^m$ has this property for all~$m$; therefore
$\Sigma^rn^m/\Sigma^rn^1$ is an $r$-reflective polynomial when $m$ is
odd, an anti-$r$-reflective polynomial when $m$ is even. In the former
case, we are done, by Lemma~4. In the latter case, Lemma~4 establishes
the existence of a polynomial~$g$ such that
$\Sigma^rn^m/\Sigma^rn^1=(n+r/2)g\bigl(n(n+r)\bigr)$. Again, we are
done, because the identity
$$\Sigma^rn^2={2n+r\over r+2}\,\Sigma^rn^1$$
is readily verified.\quad\pfbox

\bigskip\noindent
{\bf A plausible derivation.}\enspace
Faulhaber probably didn't think about $r$-reflective and
anti-$r$-reflective functions in exactly the way we have described
them, but his book~[2] certainly indicates that he was quite familiar
with the territory encompassed by that theory.

In fact, he could have found his formulas for power sums without
knowing the theory in detail. A~simple approach, illustrated here for
$\Sigma n^{13}$, would suffice: Suppose
$$14\Sigma n^{13}=n^7(n+1)^7-S(n)\,,$$
where $S(n)$ is a 1-reflective function to be determined. Then
$$\eqalign{14n^{13}&=n^7(n+1)^7-(n-1)^7n^7-\nabla S(n)\cr
\noalign{\smallskip}
&=14n^{13}+70n^{11}+42n^9+2n^7-\nabla S(n)\,,\cr}$$
and we have
$$S(n)=70\Sigma n^{11}+42n^9+2\Sigma n^7\,.$$
In other words
$$\Sigma n^{13}={64\over 7}N^7-5\Sigma n^{11}-3\Sigma n^9-{1\over
7}\Sigma n^7,$$
and we can complete the calculation by subtracting multiples of
previously computed results.

The great advantage of using polynomials in $N$ rather than $n$ is
that the new formulas are considerably shorter. The method Faulhaber
and others had used before making this discovery was most likely
equivalent to the laborious calculation
$$\eqalignno{\Sigma n^{13}&={\textstyle{1\over 14}n^{14}+{13\over
2}\Sigma n^{12}-26\Sigma n^{11}+{143\over 2}\Sigma n^{10}-143\Sigma
n^9+{429\over 2}\Sigma n^8}\cr
\noalign{\smallskip}
&\quad\null+{\textstyle{1716\over 7}\Sigma n^7+{429\over 2}\Sigma
n^6-143\Sigma n^5+{143\over 2}\Sigma n^4-26\Sigma n^3+{13\over
2}\Sigma n^2 -\Sigma n^1+{1\over 14}} n\,;\cr}$$
the coefficients here are ${1\over 14}{14\choose 12},\,-{1\over
14}{14\choose 11},\,\ldots\,,\,{1\over 14}{14\choose 0}$.

To handle sums of even exponents, Faulhaber knew that
$$\Sigma n^{2m}={n+{1\over 2}\over
2m+1}\,(a_1N+a_2N^2+\cdots+a_mN^m)$$
holds if and only if
$$\Sigma n^{2m+1}={a_1\over 2}N^2+{a_2\over 3}N^3+\cdots+{a_m\over
m+1}N^{m+1}\,.$$ 
Therefore he could get two sums for the price of one~[2, folios C.iv
verso and D.i recto]. It is not difficult to prove this relation by
establishing an isomorphism between the calculations of $\Sigma
n^{2m+1}$ and the calculations of the quantities
$\left.S_{2m}=\bigl((2m+1)\Sigma n^{2m}\bigr)\right/(n+{1\over 2})$;
for example, the recurrence for $\Sigma n^{13}$ above corresponds to
the formula
$$S_{12}=64N^6-5S_{10}-3S_8-{\textstyle{1\over 7}}S_6\,,$$
which can be derived in essentially the same way. Since the
recurrences are essentially identical, we obtain a correct formula for
$\Sigma n^{2m+1}$ from the formula for $S_{2m}$ if we replace~$N^k$ by
$N^{k+1}\!/(k+1)$.

\bigskip\noindent
{\bf Faulhaber's cryptomath.}\enspace
Mathematicians of Faulhaber's day tended to conceal their methods and
hide results in secret code. Faulhaber ends his book~[2] with a
curious exercise of this kind, evidently intended to prove to
posterity that he had in fact computed the formulas for sums of powers
as far as $\Sigma n^{25}$ although he published the results only up to
$\Sigma n^{17}$.

His puzzle can be translated into modern notation as follows: Let
$$\Sigma^9n^8={a_{17}n^{17}+\cdots+a_2n^2+a_1n\over d}$$
where the $a$'s are integers having no common factor and
$d=a_{17}+\cdots+a_2+a_1$. Let
$$\Sigma n^{25}={A_{26}n^{26}+\cdots +A_2n^2+A_1n\over D}$$
be the analogous formula for $\Sigma n^{25}$. Let
$$\eqalign{\Sigma n^{22}&={(b_{10}N^{10}-b_9N^9+\cdots +b_0)\over
b_{10}-b_9+\cdots +b_0}\;\Sigma n^2\,,\cr
\noalign{\smallskip}
\Sigma n^{23}&={(c_{10}n^{10}-c_9N^9+\cdots+c_0)\over
c_{10}-c_9+\cdots +c_0}\Sigma n^3\,,\cr
\noalign{\smallskip}
\Sigma n^{24}&={(d_{11}n^{11}-d_{10}N^{10}+\cdots-d_0)\over
d_{11}-d_{10}+\cdots +d_0}\Sigma n^2\,,\cr
\noalign{\smallskip}
\Sigma n^{25}&={(e_{11}n^{11}-e_{10}N^{10}+\cdots-e_0)\over
e_{11}-e_{10}+\cdots +e_0}\Sigma n^3\,,\cr}$$
where the integers $b_k$, $c_k$, $d_k$, $e_k$ are as small as possible
so that $b_k,c_k,d_k,e_k$ are multiples of~$2^k$. (He wants them to be
multiples of~$2^k$ so that $b_kN^k$, $c_kN^k$, $d_kN^k$, $e_kN^k$ are
polynomials in~$n$ with integer coefficients; that is why he wrote,
for example, $\Sigma n^7=(12N^2-8N+2)N^2\!/6$ instead of
$(6N^2-4N+1)N^2\!/3$. See~[2, folio D.i verso].) Then compute
$$\eqalign{x_1&=(c_3-a_{12})/7924252\,;\cr
\noalign{\smallskip}
x_2&=(b_5+a_{10})/112499648\,;\cr
\noalign{\smallskip}
x_3&=(a_{11}-b_9-c_1)/2945002\,;\cr
\noalign{\smallskip}
x_4&=(a_{14}+c_7)/120964\,;\cr
\noalign{\smallskip}
x_5&=(A_{26}a_{11}-D+a_{13}+d_{11}+e_{11})/199444\,.\cr}$$
These values $(x_1,x_2,x_3,x_4,x_5)$ specify the five letters of a
``hochger\"uhmte Nam,'' if we use five designated alphabets [2, folio
F.i recto].

It is doubtful whether anybody solved this puzzle during the first
360~years after its publication, but the task is relatively easy with
modern computers. We have
$$\vcenter{\halign{%
$#$\hfil\ &$#$\hfil\ &$#$\hfil\ &$#$\hfil\ &$#$\hfil\cr
a_{10}=532797408\,,&a_{11}=104421616\,,&a_{12}=14869764\,,%
&a_{13}=1526532\,,&a_{14}=110160\,;\cr
\noalign{\smallskip}
b_5=29700832\,,&b_9=140800\,;\cr
\noalign{\smallskip}
c_1=205083120\,,&c_3=344752128\,,&c_7=9236480\,;\cr
\noalign{\smallskip}
d_{11}=559104\,;&e_{11}=86016\,;&A_{26}=42\,;&D=1092\,.\cr}}$$
The fact that $x_2=(29700832+532797408)/112499648=5$ is an integer is
reassuring: We must be on the right track! But alas, the other values
are not integral.

A bit of experimentation soon reveals that we do obtain good results
if we divide all the~$c_k$ by~4. Then, for example,
$x_1=(344752128/4-14869764)/7924252=9$, and we also find $x_3=18$,
$x_4=20$. It appears that Faulhaber calculated $\Sigma^9n^8$ and
$\Sigma n^{22}$ correctly, and that he also had a correct expression for
$\Sigma n^{23}$ as a polynomial in~$N$; but he probably never went on
to express $\Sigma n^{23}$ as a polynomial in~$n$, because he would
then have multiplied his coefficients by~4 in order to compute
$c_6N^6$ with integer coefficients.

The values of $(x_1,x_2,x_3,x_4)$ correspond to the letters 
I$\,$E$\,$S$\,$U, so
the concealed name in Faulhaber's riddle is undoubtedly
I$\,$E$\,$S$\,$U$\,$S (Jesus).

But his formula for $x_5$ does not check out at all; it is way out of
range and not an integer.
This is the only formula that relates to $\Sigma n^{24}$ and $\Sigma
n^{25}$, and it involves only the simplest elements of those sums---the
leading coefficients $A_{26}$, $D$, $d_{11}$, $e_{11}$. Therefore we
have no evidence that Faulhaber's  calculations beyond $\Sigma n^{23}$ were
reliable. It is tempting to imagine that he meant to say
`$A_{26}a_{11}/D$' instead of `$A_{26}a_{11}-D$' in his formula
for~$x_5$, but even then major corrections are needed to the other
terms and it is unclear what he intended.

\bigskip\noindent
{\bf All-integer formulas.}\enspace
Faulhaber's theorem allows us to express the power sum $\Sigma n^m$ in
terms of about ${1\over 2}m$~coefficients. The elementary theory above
also suggests another approach that produces a similar effect: We can
write, for example,
$$\eqalign{n&={\textstyle{n\choose 1}}\,;\cr
\noalign{\smallskip}
n^3&=6{\textstyle{n+1\choose 3}+{n\choose 1}}\,;\cr
\noalign{\smallskip}
n^5&=120{\textstyle{n+2\choose 5}+30{n+1\choose 3}+{n\choose 1}}\,;\cr}$$
(It is easy to see that any odd function $g(n)$ of the integer~$n$ can
be expressed uniquely as a linear combination
$$g(n)=\textstyle{a_1{n\choose 1}
+a_3{n+1\choose 3}+a_5{n+2\choose 5}}+\cdots$$
of the odd functions ${n\choose 1},\, {n+1\choose 3},\,{n+2\choose
5}$, \dots, because we can determine the coefficients
$a_1,a_3,a_5,\ldots$ successively by plugging in the values
$n=1,2,3$, \dots. The coefficients~$a_k$ will be integers iff $g(n)$
is an integer for all~$n$.) Once $g(n)$ has been expressed in this
way, we clearly have
$$\Sigma g(n)=\textstyle{a_1{n+1\choose 2}+a_3{n+2\choose
4}+a_5{n+3\choose 6}}+\cdots\;.$$

This approach therefore yields the following identities for sums of
odd powers:
$$\vcenter{\halign{$#$\hfil\ &$#$\hfil\cr
\Sigma n^1&={n+1\choose 2}\,;\cr
\noalign{\smallskip}
\Sigma n^3&=6{n+2\choose 4}+{n+1\choose 2}\,;\cr
\noalign{\smallskip}
\Sigma n^5&=120{n+3\choose 6}+30{n+2\choose 4}+{n+1\choose 2}\,;\cr
\noalign{\smallskip}
\Sigma n^7&=5040{n+4\choose 8}+1680{n+3\choose 6}+126{n+2\choose 4}+
{n+1\choose 2}\,;\cr
\noalign{\smallskip}
\Sigma n^9&=362880{n+5\choose 10}+151200{n+4\choose 8}
+17640{n+3\choose 6}+510{n+2\choose 4}+{n+1\choose 2}\,;\cr
\noalign{\smallskip}
\Sigma n^{11}&=39916800{n+6\choose 12}+
19958400{n+5\choose 10}+3160080{n+4\choose 8}
+168960{n+3\choose 6}\cr
\noalign{\smallskip}
&\qquad\null+2046{n+2\choose 4}+{n+1\choose 2}\,;\cr
\noalign{\smallskip}
\Sigma n^{13}&=6227020800{n+7\choose 14}+3632428800{n+6\choose 12}+
726485760{n+5\choose 10}\cr
\noalign{\smallskip}
&\qquad\null+57657600{n+4\choose 8}
+1561560{n+3\choose 6}+8190{n+2\choose 4}+{n+1\choose 2}\,.\cr}}$$
And repeated sums are equally easy; we have
$$\Sigma^rn^1=\textstyle{{n+r\choose 1+r}\,,\qquad
\Sigma^rn^3=6{n+1+r\choose 3+r}+{n+r\choose 1+r}}\,,\quad\hbox{etc.}$$

The coefficients in these formulas are related to what Riordan [R,
page 213] has called {\it central factorial numbers\/}
 of the second kind. In his notation
$$x^m=\sum_{k=1}^mT(m,k)x^{[k]}\,,\quad
x^{[k]}=x\bigl(\textstyle{x+{k\over 2}
-1}\bigr)\bigl(\textstyle{x+{k\over 2}-2}\bigr)\,\ldots\,
\bigl(\textstyle{x+{k\over 2}+1}\bigr)\,,$$
when $m>0$, and $T(m,k)=0$ when $m-k$ is odd; hence
$$\eqalign{n^{2m-1}&=\sum_{k=1}^m\,(2k-1)!\,T(2m,2k){n+k-1\choose
2k-1}\,,\cr
\noalign{\smallskip}
\Sigma n^{2m-1}&=\sum_{k=1}^m\,(2k-1)!\,T(2m,2k){n+k\choose
2k}\,.\cr}$$ 
The coefficients $T(2m,2k)$ are always integers, because 
the identity $x^{[k+2]}=x^{[k]}(x^2-k^2\!/4)$ implies the recurrence
$$T(2m+2,2k)=k^2T(2m,2k)+T(2m,2k-2)\,.$$
The generating function for these numbers turns out to be
$$\hbox{cosh}\bigl(2 x \sinh(y/2)\bigr)=\sum_{m=0}^{\infty}
\biggl(\sum_{k=0}^m T(2m,2k)x^{2k}\biggr)\,{y^{2m}\over (2m)!}\,.$$

Notice that the power-sum formulas obtained in this way are more
``efficient'' than the well-known formulas based on Stirling numbers
(see~[5, (6.12)]):
$$\Sigma n^m=\sum_k k!{m\brace k}{n+1\choose k+1}=\sum_kk!{m\brace
k}(-1)^{m-k}{n+k\choose k+1}\,.$$
The latter formulas give, for example,
$$\eqalign{%
\Sigma n^7&=5040{\textstyle{n+1\choose 8}+15120{n+1\choose 7}
+16800{n+1\choose 6}+8400{n+1\choose 5}+1806{n+1\choose 4}}\cr
\noalign{\smallskip}
&\qquad\null+126{\textstyle{n+1\choose 3}+{n+1\choose 2}}\cr
\noalign{\smallskip}
&=5040{\textstyle{n+7\choose 8}-15120{n+6\choose 7}+16800{n+5\choose
6}-8400{n+4\choose 5}+1806{n+3\choose 4}}\cr
\noalign{\smallskip}
&\qquad\null-126{\textstyle{n+2\choose 3}+{n+1\choose 2}}\,.\cr}$$
There are about twice as many terms, and the coefficients are larger.
(The Faulhaberian expression $\Sigma n^7=(6N^4-4N^3+N^2)/3$ is, of course,
better yet.)

Similar formulas for even powers can be obtained as follows. We have
$$\vcenter{\halign{\hfil$#$\hfil$\;$&$#$\hfil$\;$&$#$\hfil\cr
n^2&=n{n\choose 1}&=U_1(n)\,,\cr
\noalign{\smallskip}
n^4&=6n{n+1\choose 3}+n{n\choose 1}&=12U_2(n)+U_1(n)\,,\cr
\noalign{\smallskip}
n^6&=120n{n+2\choose 5}+30n{n+1\choose 3}+n{n\choose 1}
&=360U_3(n)+60U_2(n)+U_1(n)\,,\cr}}$$
etc., where
$$U_k(n)={n\over k}\,{n+k-1\choose 2k-1}={n+k\choose
2k}+{n+k-1\choose 2k}\,.$$
Hence
$$\vcenter{\halign{$#$\hfil\ &$#$\hfil\cr
\Sigma n^2&=T_1(n)\,,\cr
\noalign{\smallskip}
\Sigma n^4&=12T_2(n)+T_1(n)\,,\cr
\noalign{\smallskip}
\Sigma n^6&=360T_3(n)+60T_2+T_1(n)\,,\cr
\noalign{\smallskip}
\Sigma n^8&=20160T_4(n)+5040T_3(n)+252T_2(n)+T_1(n)\,,\cr
\noalign{\smallskip}
\Sigma n^{10}&=1814400T_5(n)+604800T_4(n)+52920T_3(n)+1020T_2(n)+T_1(n)\,,\cr
\noalign{\smallskip}
\Sigma n^{12}&=239500800T_6(n)+99792000T_5(n)+12640320T_4(n)\cr
\noalign{\smallskip}
&\qquad\null+506880T_3(n)+4092T_2(n)+T_1(n)\,,\cr}}$$
etc., where
$$T_k(n)={n+k+1\choose 2k+1}+{n+k\choose 2k+1}={2n+1\over
2k+1}\,{n+k\choose 2k}\,.$$

Curiously, we have found a relation here between $\Sigma n^{2m}$ and
$\Sigma n^{2m-1}$, somewhat analogous to Faulhaber's relation between
$\Sigma n^{2m}$ and $\Sigma n^{2m+1}$: The formula
$${\Sigma n^{2m}\over 2n+1}={a_1{n+1\choose
2}+a_2{n+2\choose 4}+\,\cdots\,+\,a_m{n+m\choose 2m}}$$
holds if and only if
$$\Sigma n^{2m-1}={3\over 1}a_1{n+1\choose 2}+{5\over
2}a_2{n+2\choose 4}+\cdots +{2m+1\over m}a_m{n+m\choose 2m}\,.$$

\bigskip\noindent
{\bf 4. Reflective decomposition.}\enspace
The forms of the expressions in the previous section lead naturally to
useful representations of arbitrary 
functions $f(n)$ defined on the integers. It
is easy to see that any $f(n)$ can be written uniquely in the
form
$$f(n)=\sum_{k\geq0}a_k{n+\lfloor k/2\rfloor\choose k}\,,$$
for some coefficients $a_k$; indeed, we have
$$a_k=\nabla^kf(\lfloor k/2\rfloor)\,.$$
(Thus $a_0=f(0)$, $a_1=f(0)-f(-1)$, $a_2=f(1)-2f(0)+f(-1)$, etc.) The
$a_k$ are integers iff $f(n)$ is always an integer. The $a_k$ are
eventually zero iff $f$ is a polynomial. The $a_{2k}$ are all zero iff
$f$ is odd. The $a_{2k+1}$ are all zero iff $f$ is 1-reflective.

Similarly, there is a unique expansion
$$f(n)=b_0T_0(n)+b_1U_1(n)+b_2T_1(n)+b_3U_2(n)+b_4T_2(n)+\cdots\;,$$
in which the $b_k$ are integers iff $f(n)$ is always an integer. The
$b_{2k}$ are all zero iff $f$ is even and $f(0)=0$. The $b_{2k+1}$ are
all zero iff $f$ is anti-1-reflective. Using the recurrence relations
$$\nabla T_k(n)=U_k(n)\,,\qquad \nabla U_k(n)=T_{k-1}(n-1)\,,$$
we find
$$a_k=\nabla^kf(\lfloor k/2\rfloor)=2b_{k-1}+(-1)^kb_k$$
and therefore
$$b_k=\sum_{j=0}^k(-1)^{\lceil j/2\rceil+\lfloor
k/2\rfloor}2^{k-j}a_j\,.$$
In particular, when $f(n)=1$ for all~$n$, we have $b_k=(-1)^{\lfloor
k/2\rfloor} 2^k$. The infinite series is finite for each~$n$.

\proclaim
Theorem. If $f$ is any function defined on the integers and if $r,s$
are arbitrary integers, we can always express $f$ in the form
$$f(n)=g(n)+h(n)$$
where $g(n)$ is $r$-reflective and $h(n)$ is anti-$s$-reflective. This
representation is unique, except when $r$ is even and $s$ is odd; in
the latter case the representation is unique if we specify the value
of~$g$ or~$h$ at any point.

\proof
It suffices to consider $0\leq r,s\leq 1$, because $f(x)$ is
(anti)-$r$-reflective iff $f(x+a)$ is (anti)-$(r+2a)$-reflective. 

When $r=s=0$, the result is just the well known decomposition of a
function into even and odd parts,
$$g(n)={\textstyle{1\over 2}}\bigl(f(n)+f(-n)\bigr)\,,\qquad
h(n)={\textstyle{1\over 2}}\bigl(f(n)-f(-n)\bigr)\,.$$
When $r=s=1$, we have similarly
$$g(n)={\textstyle{1\over 2}}\bigl(f(n)+f(-1-n)\bigr)\,,\qquad
h(n)={\textstyle{1\over 2}}\bigl(f(n)-f(-1-n)\bigr)\,.$$

When $r=1$ and $s=0$, it is easy to deduce that $h(0)=0$, $g(0)=f(0)$,
$h(1)=f(0)-f(-1)$, $g(1)=f(1)-f(0)+f(-1)$,
$h(2)=f(1)-f(0)+f(-1)-f(-2)$, $g(2)=f(2)-f(1)+f(0)-f(-1)+f(-2)$, etc.

And when $r=0$ and $s=1$, the general solution is
$g(0)=f(0)-C$, $h(0)=C$, $g(1)=f(-1)+C$, $h(1)=f(1)-f(-1)-C$,
$g(2)=f(1)-f(-1)+f(-2)-C$, $h(2)=f(2)-f(1)+f(-1)-f(-2)+C$,
etc.\quad\pfbox

\bigskip
When $f(n)=\sum_{k\geq 0}a_k{n+\lfloor k/2\rfloor\choose k}$, the case
$r=1$ and $s=0$ corresponds to the decomposition
$$g(n)=\sum_{k=0}^{\infty}a_{2k}{n+k\choose 2k}\,,\qquad
h(n)=\sum_{k=0}^{\infty}a_{2k+1}{n+k\choose 2k+1}\,.$$
Similarly, the representation $f(n)=\sum_{k\geq
0}b_{2k}T_k(n)+\sum_{k\geq 0}b_{2k+1}U_{k+1}(n)$ corresponds to the
case $r=0$, $s=1$, $C=f(0)$.


\bigskip\noindent
{\bf Back to Faulhaber's form.}\enspace
Let us now return to representations of $\Sigma n^m$ as polynomials in
$n(n+1)$. Setting $u=2N=n^2+n$, we have
$$\vcenter{\halign{$#$\hfil\ &$#$\hfil\ &$#$\hfil\cr
\Sigma n&={1\over 2}u&={1\over 2}A_0^{(1)}u\cr
\noalign{\smallskip}
\Sigma n^3&={1\over 4}u^2&={1\over 4}\bigl(A_0^{(2)}u^2+A_1^{(2)}u\bigr)\cr
\noalign{\smallskip}
\Sigma n^5&={1\over 6}\bigl(u^3-{1\over 2}u^2\bigr)
&={1\over 6}\bigl(A_0^{(3)}u^3+A_1^{(3)}u^2+A_2^{(3)}u\bigr)\cr
\noalign{\smallskip}
\Sigma n^7&={1\over 8}\bigl(u^4-{4\over 3}u^3+{2\over 3}u^2\bigr)
&={1\over
8}\bigl(A_0^{(4)}u^4+A_1^{(4)}u^3+A_2^{(4)}u^2+A_3^{(4)}u\bigr)\cr}}$$
and so on, 
for certain coefficients $A_k^{(m)}$.

Faulhaber never discovered the Bernoulli numbers; i.e., he never
realized that a single sequence of constants $B_0,B_1,B_2,\ldots$ would
provide a uniform formula
$$\Sigma n^m={\textstyle{1\over
m+1}}\bigl(B_0n^{m+1}-{\textstyle{m+1\choose
1}}B_1n^m+{\textstyle{m+1\choose
2}}B_2n^{m-1}-\cdots+(-1)^m{\textstyle{m+1\choose m}}B_mn\bigr)$$
for all sums of powers. He never mentioned, for example, the fact that
almost half of the coefficients turned out to be zero after he had
converted his formulas for $\Sigma n^m$ from polynomials in~$N$ to
polynomials in~$n$. (He did notice that the coefficient of~$n$ was
zero when $m>1$ was odd.)

However, we know now that Bernoulli numbers exist, and we know that
$B_3=B_5=B_7=\cdots =0$. This is a strong condition. Indeed, it
completely defines the constants $A_k^{(m)}$ in the Faulhaber
polynomials above, given that $A_0^{(m)}=1$.

For example, let's consider the case $m=4$, i.e., the formula for
$\Sigma n^7$:
We need to find coefficients $a=A_1^{(4)}$, $b=A_2^{(4)}$,
$c=A_3^{(4)}$ such that the polynomial
$$n^4(n+1)^4+an^3(n+1)^3+bn^2(n+1)^2+cn(n+1)$$
has vanishing coefficients of $n^5$, $n^3$, and~$n$. The polynomial is
$$\vcenter{\halign{\hfil$#\;$&$#$\hfil$\;$&$#$\hfil$\;$&$#$\hfil\cr
n^8+4n^7&\null+6n^6+4n^5&\null+n^4\cr
&\null+an^6+3an^5&\null+3an^4+an^3\cr
&&\null+bn^4+2bn^3&\null+n^2\cr
&&&\null+cn^2+cn\,;\cr}}$$
so we must have $3a+4=2b+a=c=0$. In general the coefficient of, say,
$n^{2m-5}$ in the polynomial for $2m\Sigma n^{2m-1}$ is easily seen to
be
$$\textstyle{m\choose 5}A_0^{(m)}+{m-1\choose 3}A_1^{(m)}+{m-2\choose
1}A_2^{(m)}\,.$$
Thus the Faulhaber coefficients can be defined by the rules
$$A_0^{(w)}=1\,;\qquad\sum_{j=0}^k{w-j\choose
2k+1-2j}A_j^{(w)}=0\,,\quad k>0\,.\eqno(\ast)$$
(The upper parameter will often be called $w$ instead of $m$, in the sequel,
because we will want to generalize to noninteger values.) Notice that
$(\ast)$ defines the coefficients for each exponent without reference
to other exponents; for every integer $k\geq 0$, the quantity
$A_k^{(w)}$ is a certain rational function of~$w$. For example, we
 have
$$\eqalign{-A_1^{(w)}&=w(w-2)/6\,,\cr
\noalign{\smallskip}
A_2^{(w)}&=w(w-1)(w-3)(7w-8)/360\,,\cr
\noalign{\smallskip}
-A_3^{(w)}&=w(w-1)(w-2)(w-4)(31w^2-89w+48)/15120\,,\cr
\noalign{\smallskip}
A_4^{(w)}&=w(w-1)(w-2)(w-3)(w-5)(127w^3-691w^2+1038w-384)/6048000\,,\cr}$$
and in general $A_k^{(w)}$ is
$w^{\underline{k}}=w(w-1)\,\ldots\,(w-k+1)$
 times a polynomial of degree~$k$, with leading coefficient equal to
 $(2-2^{2k})B_{2k}/(2k)!\,$; if
$k>0$, that polynomial vanishes when $w=k+1$. 

Jacobi mentioned these coefficients $A_k^{(m)}$ in~[6],
although he did not consider the recurrence~($\ast$),
and he tabulated them for $m\leq 6$. He observed that the derivative of
$\Sigma n^m$ with respect to~$n$ is $m\,\Sigma n^{m-1}+B_m$; this
follows because power sums can be expressed in terms of Bernoulli
polynomials, 
$$\Sigma n^m={\textstyle{1\over m+1}}\bigl(B_{m+1}(n+1)-B_{m+1}(0)\bigr)\,,$$
and because $B'_m(x)=mB_{m-1}(x)$. Thus Jacobi obtained a new proof of
Faulhaber's formulas for even exponents
$$\eqalign{\Sigma n^2&={\textstyle{1\over 3}\bigl({2\over
4}A_0^{(2)}u+{1\over 4}A_1^{(2)}\bigr)(2n+1)\,,}\cr
\noalign{\smallskip}
\Sigma n^4&={\textstyle{1\over 5}\bigl({3\over 6}A_0^{(3)}u^2+{2\over
6}A_1^{(3)}u+{1\over 6}A_2^{(3)}\bigr)(2n+1)\,,}\cr
\noalign{\smallskip}
\Sigma n^6&={\textstyle{1\over 7}\bigl({4\over 8}A_0^{(4)}u^3+{3\over
8}A_1^{(4)}u^2 +{2\over 8}A_2^{(4)}u+{1\over
8}A_3^{(4)}\bigr)(2n+1)\,,}\cr}$$ 
etc. (The constant terms are zero, but they are shown explicitly here
 so that the pattern is plain.) Differentiating again gives, e.g.,
$$\eqalign{\Sigma n^5&={1\over 6 \cdot 7\cdot 8}
\bigl((4\cdot 3\,A_0^{(4)}u^2+3\cdot 2\,
A_1^{(4)}u+2\cdot 1\,A_2^{(4)})(2n+1)^2\cr
\noalign{\smallskip}
&\qquad\qquad
\qquad\null+2(4A_0^{(4)}u^3+3A_1^{(4)}u^2+2A_2^{(4)}u+1A_3^{(4)})\bigr)-
{\textstyle{1\over 6}}B_6\cr
\noalign{\smallskip}
&={1\over 6\cdot 7\cdot 8}\bigl(8\cdot 7\,A_0^{(4)}u^3+(6\cdot
5\,A_1^{(4)}+4\cdot 3\,A_0^{(4)})u^2+(4\cdot 3\,A_2^{(4)}+3\cdot 2\,
A_1^{(4)})u\cr
\noalign{\smallskip}
&\qquad\qquad\qquad\null+(2\cdot 1\,A_3^{(4)}+2\cdot
1\,A_2^{(4)})\bigr)-{\textstyle{1\over 6}}B_6\,.\cr}$$
This yields Jacobi's recurrence
$$(2w-2k)(2w-2k-1)A_k^{(w)}+(w-k+1)(w-k)A_{k-1}^{(w)}
=2w(2w-1)A_k^{(w-1)}\,,\eqno(\ast\ast)$$
which is valid for all integers $w>k+1$ so it must be valid for
all~$w$. Our derivation of ($\ast\ast$) also allows us to  conclude that
$$A_{m-2}^{(m)}={2m\choose 2}B_{2m-2}\,,\quad m\geq 2\,,$$
by considering the constant term of the second derivative of $\Sigma
n^{2m-1}$. 

Recurrence $(\ast)$ does not define $A_m^{(m)}$, except as the limit
of $A_m^{(w)}$ when $w\rightarrow m$. But we can compute this value
by setting $w=m+1$ and $k=m$ in~$(\ast\ast)$, which reduces to
$$2A_{m-1}^{(m+1)}=(2m+2)(2m+1)A_m^{(m)}$$
because $A_m^{(m+1)}=0$. Thus
$$A_m^{(m)}=B_{2m}\,,\quad \hbox{\rm integer }m\geq 0\,.$$

\bigskip\noindent
{\bf Solution to the recurrence.}\enspace
	An explicit formula for $A_k^{(m)}$ can be found as follows: We have
$$\Sigma n^{2m-1}={1\over 2m}\bigl(B_{2m}(n+1)-B_{2m}\bigr)={1\over
2m} (A_0^{(m)}u^m+\cdots +A_{m-1}^{(m)}u)\,,$$
and $n+1=(\sqrt{\mathstrut 1+4u}+1)/2$; hence, using the known values
of~$A_m^{(m)}$, we obtain
$$\sum_{k=0}^{\infty}A_k^{(m)}u^{m-k}=B_{2m}
\left({\sqrt{\mathstrut 1+4u}+1\over
2}\right)=B_{2m}\left({1-\sqrt{\mathstrut 1+4u}\over 2}\right)\,,$$
a closed form in terms of Bernoulli polynomials.
$\bigl($We have used the fact that $A_{m+1}^{(m)}=A_{m+2}^{(m)}=\cdots=0$,
together with the identity
%
$B_n(x+1)=(-1)^nB_n(-x)\,.\bigr)$
Expanding the right side in powers of $u$ gives
$$\eqalign{&\sum_l{2m\choose l}\left({1-\sqrt{\mathstrut 1+4u}\over
2}\right)^lB_{2m-l}\cr
\noalign{\smallskip}
&\qquad\qquad =\sum_{j,l}{2m\choose l}{2j+l\choose j}{l\over
2j+l}(-u)^{j+l}B_{2m-l}\,,\cr}$$
using equation (5.70) of~[5]. Setting $j+l=m-k$ finally yields
$$A_k^{(m)}=(-1)^{m-k}\sum_j{2m\choose m-k-j}{m-k+j\choose
j}{m-k-j\over m-k+j}B_{m+k+j}\,,\quad 0\leq k<m\,.$$
This formula, which was first obtained by Gessel and Viennot~[4],
makes it easy to confirm that $A_{m-1}^{(m)}=0$ and
$A_{m-2}^{(m)}={2m\choose 2}B_{2m-2}$, and to derive additional values
such as
$$\eqalign{A_{m-3}^{(m)}&=-2{2m\choose
2}B_{2m-2}=-2A_{m-2}^{(m)}\,,\quad m\geq 3\,;\cr
\noalign{\smallskip}
A_{m-4}^{(m)}&={2m\choose 4}B_{2m-4}+5{2m\choose 2}B_{2m-2}\,,\quad
m\geq 4\,.\cr}$$

The author's interest in Faulhaber polynomials was inspired by the
work of Edwards~[1], who resurrected Faulhaber's work after it had
been long forgotten and undervalued by historians of mathematics.
Ira Gessel responded
to the same stimulus by submitting problem E3204 to the {\sl Math
Monthly\/}~[3]
regarding a bivariate generating function for
Faulhaber's coefficients. Such a function is obtainable from the
univariate generating function above, using the standard generating
function for Bernoulli polynomials: Since
$$\eqalign{\sum B_{2m}\left({x+1\over 2}\right){z^{2m}\over (2m)!}
&={1\over 2}\sum B_m\left({x+1\over 2}\right){z^m\over m!}+{1\over
2}\sum B_m\left({x+1\over 2}\right){(-z)^m\over m!}\cr
\noalign{\smallskip}
&={z\,e^{(x+1)z/2}\over 2(e^z-1)}-{z\,e^{-(x+1)z/2}\over
2(e^{-z}-1)}={z\;{\rm
cosh}(xz/2)\over 2\;{\rm sinh}(z/2)}\,,\cr}$$
we have
$$\eqalign{\sum_{k,n}A_k^{(m)}u^{m-k}\,{z^{2m}\over (2m)!}
&=\sum_mB_{2m}\left({\sqrt{\mathstrut 1+4u}+1\over 2}\right)\,
{z^{2m}\over (2m)!}
= {z\over 2}\,{{\rm cosh}\,(\sqrt{\mathstrut 1+4u}\,z/2)\over {\rm
sinh}\,(z/2)}\,;\cr
\noalign{\smallskip}
\sum_{k,w}A_k^{(m)}u^k\,{z^{2m}\over (2m)!}
&={z\sqrt{\mathstrut u}\,{\rm cosh}\,(\sqrt{\mathstrut u+4}\,z/2)\over 2\,{\rm
sinh}\,(z\sqrt{\mathstrut u}\,/2)}\,.\cr}$$

The numbers $A_k^{(m)}$ are obtainable by inverting a lower triangular
matrix, as Edwards showed; indeed, recurrence $(\ast)$ defines such a
matrix. Gessel and Viennot~[4]
observed that we can therefore express them in terms of a $k\times k$
determinant, 
$$
A_k^{(w)}={1\over (1-w)\,\ldots\,(k-w)}\,
\left\vert\vcenter{\halign{$\hfil#\hfil$\quad&$\hfil#\hfil$\quad%
&$\hfil#\hfil$\quad&$\hfil#\hfil$\quad&$\hfil#\hfil$\cr
{w-k+1\choose 3}&{w-k+1\choose 1}&0&\ldots&0\cr
\noalign{\smallskip}
{w-k+2\choose 5}&{w-k+2\choose 3}&{w-k+2\choose 1}&\ldots&0\cr
\noalign{\smallskip}
\vdots&\vdots&\vdots\cr
\noalign{\smallskip}
{w-1\choose 2k-1}&{w-1\choose 2k-3}&{w-1\choose
2k-5}&\ldots&{w-1\choose 1}\cr
\noalign{\smallskip}
{w\choose 2k+1}&{w\choose 2k-1}&{w\choose 2k-3}&\ldots&{w\choose
3}\cr}}\right\vert\,.$$
When $w$ and $k$ are positive integers, Gessel and Viennot proved that
this determinant is the number of sequences of positive integers
$a_1a_2a_3\ldots a_{3k}$ such that
$$\displaylines{a_{3j-2}<a_{3j-1}<a_{3j}\leq w-k+j\,,\quad{\rm
for}\quad 1\leq j\leq k\,,\cr
\noalign{\smallskip}
a_{3j-2}<a_{3j+1}\,,\quad a_{3j-1}<a_{3j+3}\,,\quad{\rm
for}\quad 1\leq j< k\,.\cr}$$
In other words, it is the number of ways to put positive integers into
a $k$-rowed triple staircase such as
$$\vbox{\offinterlineskip
\def\hb{\phantom{\hbox{A}}}
\halign{\strut#&\vrule#&\hfil#\hfil%
&\vrule#&\hfil#\hfil%
&\vrule#&\hfil#\hfil%
&\vrule#&\hfil#\hfil%
&\vrule#&\hfil#\hfil%
&\vrule#&\hfil#\hfil%
&\vrule#&\hfil#\hfil%
&\vrule#\cr
\omit&\omit&\omit&\omit&\omit&\omit&\omit&\omit
&\multispan7{\kern-.4pt\hrulefill\kern-.4pt}\cr
\omit&\omit&\omit&\omit&\omit&\omit&\omit&\omit&&\hb&&\hb&&\hb&&\cr
\omit&\omit&\omit&\omit&\omit&\omit
&\multispan9{\kern-.4pt\hrulefill\kern-.4pt}\cr
\omit&\omit&\omit&\omit&\omit&\omit&&\hb&&\hb&&\hb&&\cr
\omit&\omit&\omit&\omit&\multispan9{\hrulefill}\cr
\omit&\omit&\omit&\omit&\omit&\hb&&\hb&&\hb&&\cr
\omit&\omit&\multispan9{\kern-.4pt\hrulefill\kern-.4pt}\cr
\omit&\omit&\omit&\hb&&\hb&&\hb&&\cr
\multispan9{\kern-.4pt\hrulefill\kern-.4pt}\cr
\omit&\hb&&\hb&&\hb&&\cr
\multispan7{\kern-.4pt\hrulefill\kern-.4pt}\cr
}}$$
with all rows and all columns strictly increasing from left to right
and from top to bottom, and with all entries in row~$j$ at most
$w-k+j$. This provides a surprising combinatorial interpretation of
the Bernoulli number~$B_{2m}$ when $w=m+1$ and $k=m-1$ (in which case
the top row of the staircase is forced to contain $1,2,3$).

The combinatorial interpretation proves in particular that
$(-1)^kA_k^{(m)}\geq 0$ for all $k\geq 0$. Faulhaber stated this, but
he may not have known how to prove~it.

Denoting the determinant by $D(w,k)$, Jacobi's recurrence $(\ast\ast)$
implies that we have
$$\eqalign{&(w-k)^2(w-k+1)(w-k-1)D(w,k-1)\cr
&\qquad =(2w-2k)(2w-2k-1)(w-k-1)D(w,k)-2w(2w-1)(w-1)D(w-1,k)\,;\cr}$$
this can also be written in a slightly tidier form, using a special
case of the ``integer basis'' polynomials discussed above:
$$D(w,k-1)=T_1(w-k-1)D(w,k)-T_1(w-1)D(w-1,k)\,.$$
It does not appear obvious that the determinant satisfies such a
recurrence, nor that the solution to the recurrence should have
integer values when $w$ and~$k$ are integers. But, identities are not
always obvious.

\bigskip\noindent
{\bf Generalization to noninteger powers.}\enspace
Recurrence $(\ast)$ does not require $w$ to be a positive integer, and
we can in fact solve it in closed form when $w=3/2$:
$$\eqalign{\sum_{k\geq 0}A_k^{(3/2)}u^{3/2-k}
&=B_3\left({\sqrt{\mathstrut 1+4u}+1\over 2}\right)\cr
\noalign{\smallskip}
&={u\over 2}\,\sqrt{\mathstrut 1+4u}=u^{3/2}\sum_{k\geq 0}{1/2\choose
k}(4u)^{-k}\,.\cr}$$ 
Therefore $A_k^{(3/2)}={1/2\choose k}4^{-k}$ is related to the $k$th
Catalan number.
A similar closed form exists for $A_k^{(m+1/2)}$ when $m$
is any nonnegative integer.

For other cases of $w$, our generating function for $A_k^{(w)}$
involves~$B_n(x)$ with noninteger subscripts. The Bernoulli polynomials
 can be generalized to a family of functions $B_z(x)$, for
arbitrary~$z$, in several ways; the best generalization for our
present purposes seems to arise when we define
$$B_z(x)=x^z\sum_{k\geq 0}{z\choose k}x^{-k}B_k\,,$$
choosing a suitable branch of the function $x^z$. With this definition
we can develop
the right-hand side of
$$\eqalignno{\sum_{k\geq 0}A_k^{(w)}u^{-k}
&=B_{2w}\left({\sqrt{\mathstrut 1+4u}+1\over 2}\right)u^{-w}\cr
\noalign{\smallskip}
&=\left({\sqrt{\mathstrut 1+4u}+1\over 
2\sqrt{\mathstrut u}}\right)^{2w}\sum_{k\geq
0}{2w\choose k}\left({\sqrt{\mathstrut 1+4u}+1\over
2}\right)^{-k}B_k&(\ast{\ast}\ast)\cr}$$ 
as a power series in $u^{-1}$ as $u\rightarrow \infty$. 

The factor
outside the $\sum$ sign is rather nice; we have
$$\left({\sqrt{\mathstrut 1+4u}+1\over 2\sqrt{\mathstrut
u}}\right)^{2w} =\sum_{j\geq
0}\,{w\over w+j/2}\,{w+j/2\choose j}u^{-j/2}\,,$$
because the generalized binomial series $B_{1/2}(u^{-1/2})$ 
[5,~equation (5.58)]
is the solution to
$$f(u)^{1/2}-f(u)^{-1/2}=u^{-1/2}\,,\qquad f(\infty)=1\,,$$
namely
$$f(u)=\left({\sqrt{\mathstrut 1+4u}+1\over 2\sqrt{\mathstrut u}}\right)^2\,.$$
Similarly we find
$$\eqalign{\left({\sqrt{\mathstrut 1+4u}+1\over 2}\right)^{-k}
&=\sum_j\,{-k\over j-k}\,{j/2-k/2\choose j}u^{-k/2-j/2}\cr
\noalign{\smallskip}
&=u^{-k/2}-\sum_{j\geq 1}\,{k\over 2j}{j/2-k/2-1\choose
j-1}u^{-k/2-j/2}\,.\cr}$$ 
So we can indeed expand the right-hand side as a power series with
coefficients that are polynomials in~$w$. It is actually a power
series in~$u^{-1/2}$, not~$u$; but since the coefficients of odd
powers of $u^{-1/2}$ vanish when $w$ is a positive integer, they must
be identically zero. Sure enough, a~check with computer algebra on
formal power series yields
$1+A_1^{(w)}u^{-1}+A_2^{(w)}u^{-2}+A_3^{(w)}u^{-3}+O(u^{-4})$, where
the values of $A_k^{(w)}$ for $k\leq 3$ agree perfectly
with those obtained
directly from~$(\ast)$. Therefore this approach  allows
us to express $A_k^{(w)}$ as a polynomial in~$w$, using ordinary
Bernoulli number coefficients:
$$\eqalign{A_k^{(w)}&=\sum_{l=0}^{2k}\,{w\over w+l/2}{w+l/2\choose
l}\times\cr
\noalign{\smallskip}
&\qquad\qquad\biggl({2w\choose 2k-l}B_{2k-l}-{1\over 2}\,
\sum_{j=1}^{2k-l-1}{2w\choose j}\,{j\over 2k-l-j}{k-l/2-j-1\choose
2k-l-j-1}B_j\biggr)\,.\cr}$$

The power series $(\ast{\ast}\ast)$ we have used in this successful
derivation is actually divergent for all~$u$ unless $2w$ is a
nonnegative integer, because $B_k$ grows superexponentially while the
factor 
$${2w\choose k}=(-1)^k{k-2w-1\choose k}={(-1)^k\,\Gamma(k-2w)\over
\Gamma(k+1)\,\Gamma(-2w)}\,\sim\,{(-1)^k\over\Gamma(-2w)}k^{-2w-1}$$
does not decrease very rapidly as $k\rightarrow\infty$. 

Still, $(\ast{\ast}\ast)$ is easily seen to be a valid asymptotic series
as $u\rightarrow\infty$, because asymptotic series multiply like
formal power series. This means that, for any positive integer~$p$,
we have
$$\sum_{k=0}^{2p}{2w\choose k}\left({\sqrt{\mathstrut 1+4u}+1\over
2}\right)^{2w-k}\!\!B_k 
\;=\;\sum_{k=0}^pA_k^{(w)}u^{w-k}+O(u^{w-p-1})\,.$$

We can now apply these results to obtain sums of noninteger powers, as
asymptotic series of Faulhaber's type. Suppose, for example, that we
are interested in the sum
$$H_n^{(1/3)}=\sum_{k=1}^n\,{1\over k^{1/3}}\,.$$
Euler's summation formula
[5, exercise 9.27]
tells us that 
$$\eqalign{H_n^{(1/3)}-\zeta({\textstyle{1\over 3}})
&\sim {\textstyle{3\over 2}}n^{2/3}+{\textstyle{1\over
2}}n^{-1/3}-{\textstyle{1\over 36}}n^{-4/3}-\cdots\cr
\noalign{\smallskip}
&={3\over 2}\biggl(\sum_{k\geq 0}{2/3\choose
k}n^{2/3-k}B_k+n^{-1/3}\biggr)\,,\cr}$$
and the parenthesized quantity is what we have called $B_{2/3}(n+1)$.
And when $u=n^2+n$ we have
$B_{2/3}(n+1)=B_{2/3}\bigl((\sqrt{\mathstrut 1+4u}+1)/2\bigr)$; hence
$$\eqalign{H_n^{(1/3)}-\zeta({\textstyle{1\over 3}})
&\sim {\textstyle{3\over 2}}\sum_{k\geq 0}A_k^{(1/3)}u^{1/3-k}\cr
\noalign{\smallskip}
&={\textstyle{3\over 2}}\,u^{1/3}+{\textstyle{5\over
36}}\,u^{-2/3}-{\textstyle{17\over
1215}}\,u^{-5/3}+\cdots \cr}$$
as $n\rightarrow\infty$. (We can't claim that this series converges
twice as fast as the usual one, because both series diverge! But we do get
twice as much precision in a fixed number of terms.)

In general, the same argument establishes the asymptotic series
$$\sum_{k=1}^nk^{\alpha}-\zeta(-\alpha)\,\sim\,{1\over \alpha+1}\,
\sum_{k\geq
0}A_k^{((\alpha+1)/2)}\,u^{(\alpha +1)/2-k}\,,$$
whenever $\alpha\neq -1$. The series on the right is finite when
$\alpha$ is a positive odd integer; it is convergent (for sufficiently
large~$n$) if and only if $\alpha$ is a nonnegative integer.

The special case $\alpha=-2$ has historic interest, so it deserves a
special look:
$$\eqalign{\sum_{k=1}^n\,{1\over k^2}
&\sim {\pi^2\over
6}-A_0^{(-1/2)}\,u^{-1/2}
-A_1^{(-1/2)}\,u^{-3/2}-\cdots \cr
\noalign{\smallskip}
&={\pi^2\over 6}-u^{-1/2}+{5\over 24}\,u^{-3/2}
-{161\over 1920}\,u^{-5/2}
+{401\over 7168}\,u^{-7/2}
-{32021\over 491520}\,u^{-9/2}+\cdots\;.\cr}$$
These coefficients do not seem to have a simple closed form;
 the prime factorization 
$32021=11\cdot 41\cdot 71$ is no doubt just a quirky coincidence.

\bigskip\noindent
{\bf Acknowledgments.}\enspace
This paper could not have been written without the help provided by
several correspondents.
Anthony Edwards kindly sent me a photocopy of Faulhaber's {\sl
Academia Algebr{\ae}}, a~book that is evidently extremely rare: An
extensive search of printed indexes and electronic indexes indicates
that no copies have ever been recorded to exist in America, in the
British Library, or the Biblioth\`eque Nationale. Edwards found it at
Cambridge University Library, where the volume once owned by Jacobi now
resides. (I~have annotated the photocopy and deposited it in the
Mathematical Sciences Library at Stanford, so that other interested
scholars can take a look.) Ivo Schneider, who is currently preparing a
book about Faulhaber and his work, helped me understand some of the
archaic German phrases. Herb Wilf gave me a vital insight by discovering the
first half of Lemma~4, in the case $r=1$. And Ira Gessel pointed out
that the coefficients in the expansion $n^{2m+1}=\sum a_k{n+k\choose
2k+1}$ are central factorial numbers in slight disguise.

\bigskip
\centerline{\bf References}

\bib
[1]
A.~W.~F. Edwards, ``A~quick route to sums of powers,''
{\sl American Mathematical Monthly\/ \bf 93} (1986), 451--455.

\bib
[2] 
Johann Faulhaber,
{\sl Academia Algebr{\ae}}, Darinnen die miraculosische Inventiones zu
den h\"ochsten Cossen weiters {\it continuirt\/} und {\it
profitiert\/} werden. Augspurg, bey Johann Ulrich Sch\"onigs, 1631.
(Call number 
{\tt QA154.8{\kern3pt}F3{\kern3pt}1631a{\kern3pt}f{\kern3pt}MATH} 
at Stanford University Libraries.)

\bib
[3]
Ira Gessel and University of South Alabama Problem Group,
``A~formula for power sums,'' {\sl American Mathematical Monthly\/ \bf
95} (1988), 961--962.

\bib
[4]
Ira M. Gessel and G\'erard Viennot, ``Determinants, paths, and plane
partitions,'' pre\-print, 1989.

\bib
[5]
Ronald L. Graham, Donald E. Knuth, and Oren Patashnik,
{\sl Concrete Mathematics\/}
 (Reading, Mass.: Addison-Wesley, 1989).

\bib
[6] 
C.~G.~J. Jacobi, ``De usu legitimo formulae summatoriae
Maclaurinianae,'' {\sl Journal f\"ur die reine und angewandte
Mathematik\/ \bf 12} (1834), 263--272.

\bib
[7]
John Riordan, {\sl Combinatorial Identities\/} (New York: John Wiley
\& Sons, 1968).

\bib
[8]
L.~Tits, ``Sur la sommation des puissances num\'eriques,'' {\sl
Mathesis\/ \bf 37} (1923), 353--355.

\bye